\documentclass[twoside] {siamltex}
\usepackage{graphicx}
\usepackage{cite}
\usepackage{amsmath}
\usepackage[figuresright]{rotating}
\title{Data assimilation in 2D nonlinear coupled sound and heat flow, using a stabilized explicit finite difference scheme marched backward in time.}
\author{Alfred S. Carasso\thanks
{Applied and Computational Mathematics Division,
National Institute of Standards and Technology,
Gaithersburg, MD 20899. (alfred.carasso@nist.gov).}}
\begin{document}
\maketitle
\begin{abstract}
	This paper considers the ill-posed data assimilation problem associated with hyperbolic/parabolic systems describing 2D coupled sound and heat flow.
	Given hypothetical data at time $T > 0$, that {\em may not} correspond to an actual solution of the dissipative system at time $T$, initial data at time $t=0$
	are sought that can evolve,  through the dissipative system, into a useful approximation to the desired data  at time $T$. That {\em may not}
	always be possible. A stabilized {\em explicit} finite difference scheme, marching backward in time, is developed and applied to nonlinear examples in non rectangular
	regions.
Stabilization is achieved by applying  a compensating smoothing operator at each time step, to quench the instability.
Analysis of convergence is restricted to the transparent case of linear, autonomous, selfadjoint spatial differential operators.
However, the actual computational scheme
can be applied to more general problems.
Data assimilation is illustrated using $512 \times 512$ pixel images. Such images are associated with highly irregular non smooth intensity data that
severely challenge ill-posed reconstruction procedures.
Successful and unsuccessful examples are presented.



\end{abstract}
\section{Introduction}
As was the case in \cite{Burg,Advec,Thermoassim, LeapNavAssim}, the present paper considers the {\em data assimilation} problem of recreating plausible initial values at $t=0$,
given {\em hypothetical} and/or {\em partially known} data at some later time $T >0$, in 
a coupled hyperbolic/parabolic system involving ill-posed time-reversed 2D coupled sound and heat flow, \cite{richtmyer, lattes, lions, carSound}. A
particularly advantageous direct, non iterative, {\em explicit}, backward marching finite difference scheme, is constructed and explored.
There is considerable interest in data assimilation in the geophysical sciences
\cite{cintra, howard, blum, qizhi, arcucci, Antil, Chong, Lund,auroux1,ou,auroux2,auroux3,pozo,gosse,gomez, xu, tomislava,camposvelho}, where such problems are most often treated using iterative algorithms that may include neural networks coupled with machine learning. The direct methods discussed here may
provide useful initial solutions that might be further refined by such iterative procedures. Additionally, similar direct methods may be developed in specific geophysical contexts, and these can be used, when needed,  to provide confirmation of unexpected results 
obtained by artificial intelligence methods.

Here, as was emphasized in \cite{Burg,Advec,Thermoassim, LeapNavAssim}, the given hypothetical data at time $T > 0$ may
not be smooth,
{\em may not correspond} to an actual solution at time $T$, and may differ from such
a solution by an {\em unknown large $\delta > 0$} in an  appropriate ${\cal{L}}^p$ norm. Moreover, it {\em may not be possible} to
locate initial values that can evolve into a useful approximation to the desired data at time $T$. The above data assimilation problem differs fundamentally from the ill-posed
{\em backward recovery} problems discussed in \cite{carGEM, car1IPSE, car2IPSE,car3IPSE,Thermo, car5IPSE,car6IPSE,car7IPSE}, where the given data at time $T > 0$ are  noisy, but relatively smooth, and are known to approximate an {\em actual solution} at time $T$, to within a {\em known small $\delta > 0$}, in an  appropriate ${\cal{L}}^p$ norm. For ill-posed initial value problems, all consistent stepwise marching schemes, whether explicit or implicit, are necessarily
unconditionally unstable and lead to explosive error growth, \cite[p.~59]{richtmyer}. Nevertheless, it is possible to stabilize such schemes by applying an appropriate compensating smoothing operator
at each time step to quench the instability. In \cite{carGEM, car1IPSE, car2IPSE,car3IPSE,Thermo, car5IPSE,car6IPSE,car7IPSE}, such stabilized backward marching explicit schemes have been successfully applied in backward recovery problems, and on time intervals $[0,T]$ that are significantly larger than might be expected, 
based on the uncertainty estimates in \cite{karen, payne3, knopslog,knops,payne, car00,car77,hao}. As will be seen in Section 2.1 below, the data assimilation problem presents additional difficulties.  However, limited success is still feasible in that problem, using backward marching stabilized explicit schemes.

A particularly effective vehicle for computational exploration of the proposed direct explicit data assimilation approach, lies in the use of 8 bit grey scale
$512\times 512$ pixel images, as hypothetical data at time $T > 0$. As shown in Figure 1, many natural images are defined by highly non smooth intensity data
that severely challenge ill-posed reconstructions. In the experiments to be described below, three such non smooth images are involved, and they interact with each other as the evolution progresses.

\begin{figure}
        \centerline{\includegraphics[width=5.25in]{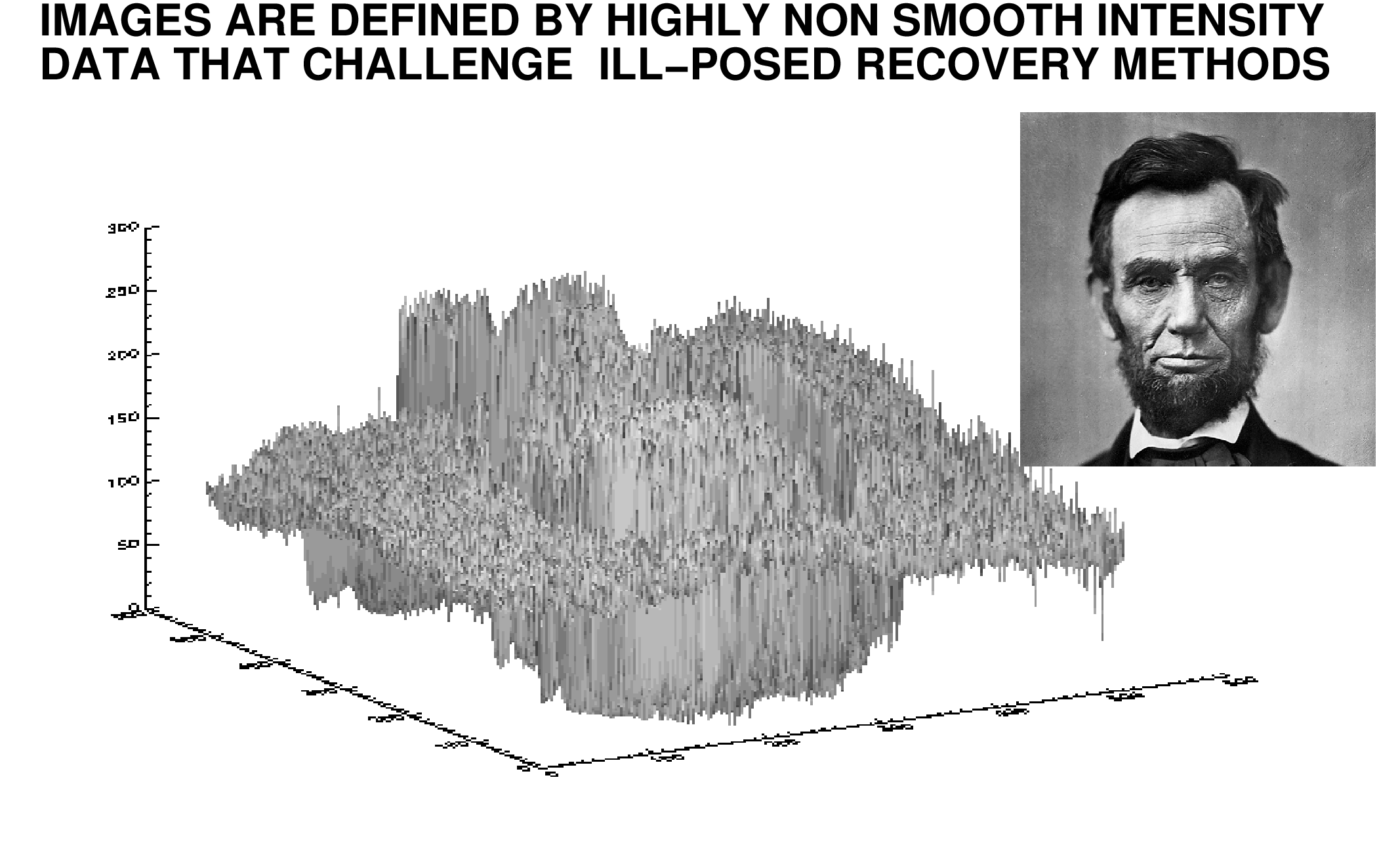}}
	\caption{FIGURE 1. Non smooth intensity data plot, associated with Abraham Lincoln image, is typical of many natural images.}
\end{figure}

\section{A stabilized explicit scheme for linear selfadjoint time-reversed coupled wave and diffusion equations}

Let $\Omega$ be a bounded domain in $R^n$ with a smooth boundary $\partial \Omega$. Let $<~,~>$ and $\parallel~\parallel_2$, respectively denote the scalar product and
norm on  ${\cal{L}}^2(\Omega)$.
Let $L$ denote a linear, second order, {\em positive definite} selfadjoint variable coefficient elliptic differential operator in $\Omega$, with homogeneous Dirichlet boundary conditions on $\partial \Omega$.
Let $\{\phi_m\}_{m=1}^{\infty}$ be the complete set of orthonormal eigenfunctions for $L$ on $\Omega$,
and let $\{\lambda_m\}_{m=1}^{\infty}$, satisfying
\begin{equation}
0 < \lambda_1 \leq \lambda_2 \leq \cdots \leq \lambda_m \leq \cdots \uparrow \infty,
\label{eq:9.00991}
\end{equation}
be the corresponding eigenvalues.

With positive constants $a,~b,~d,~$ and $0 < t \leq T_{max}$, consider the linear initial value problem on $\Omega \times (0,T_{max}]$,
\begin{equation}
\begin{array}{l}
u_t= -bLu -dv, \\
\ \\
v_t = aLu - aLw, \\
\ \\
w_t=v, \\
\ \\
u(x,0)=f(x),~~v(x,0)=g(x),~~w(x,0)=h(x).
\end{array}
\label{eq:9.0001}
\end{equation}
When  $L=-\Delta,~a=c^2,~b=\sigma,~d=(\gamma-1),~$ the above system  reduces to
the linearized equations of coupled sound and heat flow discussed in \cite{richtmyer, lattes, lions, carSound}, namely, $w_{tt}=c^2\Delta w -c^2 \Delta u,~~u_t=\sigma \Delta u -(\gamma -1) w_t~$,
with $w=u=0$ on $\partial \Omega$,
where $c$ is the isothermal sound speed, $\sigma$ is the thermal conductivity, and $~1 < \gamma < 2,~$ is the ratio of specific heats.

The initial value problem Eq.~(\ref{eq:9.0001}) becomes ill-posed when the time direction is reversed.
We contemplate such time-reversed computations by allowing for possible {\em negative} time steps $\Delta t$ in the explicit difference scheme Eq.(\ref{eq:9.0006a}) below.
With $\lambda_m$ as in Eq.~(\ref{eq:9.00991}), the positive constants $a,~b,~d,$ and the operator $L$ as in Eq.~(\ref{eq:9.0001}), ~fix $~\omega > 0$ and $p > 1$.
Given $\Delta t$, define $\rho, ~\Lambda,~Q,~\zeta_m,~q_m$, as follows:
\begin{equation}
\begin{array}{l}
\rho=\{1+d +d^2+2a^2 +2b +\sqrt{2a^2+2b^2}\},~~~\Lambda=\rho (I+ L),~~~Q=\exp(-\omega |\Delta t| \Lambda^p),\\
\ \\
\zeta_m=\rho (1+\lambda_m) > 1,~~~~~~q_m=\exp \left(-\omega |\Delta t| (\zeta_m)^p\right), \qquad m\geq 1.
\end{array}
\label{eq:9.0005a}
\end{equation}
Let $G$, $S$, and $P$, be the following $3\times 3$ matrices
\begin{equation}
G =
\begin{bmatrix}
 -bL & -dI & 0 \\
 aL & 0 & -aL \\
 0 & I & 0
\end{bmatrix}
,
\qquad
S =
\begin{bmatrix}
 Q & 0 & 0\\
 0 & Q & 0 \\
 0 & 0 & Q
\end{bmatrix}
,
\qquad
P=
\begin{bmatrix}
 \Lambda^p &  0 &  0 \\
 0 & \Lambda^p & 0  \\
 0 & 0 & \Lambda^p
\end{bmatrix}
.
\label{eq:9.0006}
\end{equation}

Let $W$ be the three component vector $[u, v, w]^T$. We may rewrite Eq.~(\ref{eq:9.0001}) as the equivalent first order system,
\begin{equation}
W_t=GW, ~~0 < t \leq T_{max}, \qquad W(\cdot, 0)=[f, g, h]^T.
\label{eq:9.0004}
\end{equation}

It is instructive to study the following explicit time-marching finite difference scheme for Eq.(\ref{eq:9.0004}), in which only the time variable is discretized, while
the space variables remain continuous. With a given positive integer $N$, let $|\Delta t|=T_{max}/N$ be the time step magnitude, and let $W^n$ denote $W(\cdot, n \Delta t), ~n=0,1,\cdots N$.
If $W(\cdot, t)$ is the unique solution of Eq.(\ref{eq:9.0004}), then
\begin{equation}
W^{n+1}=W^n + \Delta t GW^n + \tau^n,
\label{eq:9.0004a}
\end{equation}
where the `truncation error' $\tau^n = \frac{1}{2}(\Delta t)^2 G^2W(\tilde{t})$, with $n |\Delta t|  < \tilde {t} < (n+1) |\Delta t|$.
With $G$ and $S$ as in Eq.(\ref{eq:9.0006}), let $R$ be the linear operator $R = S+\Delta t SG$. We consider approximating $W^n$ with $U^n \equiv [u^n,v^n]^T$, where
\begin{equation}
U^{n+1}=S U^n + \Delta t S G U^n \equiv RU^n,~~n=0,1, \cdots (N-1), \qquad U^0=[f,g, h]^T.
\label{eq:9.0006a}
\end{equation}
\ \\
{\bf Remark.}
While the analysis that follows assumes $L$ to be a linear operator, the stabilized explicit scheme can actually be used with {\em nonlinear} operators $L$, by lagging the nonlinearity at the previous time step.
This is the case in the computational experiments to be discussed below. \\
\ \\
With $\Delta t >0$ and the data $U^0$ at time $t=0$, the forward marching scheme in Eq.(\ref{eq:9.0006a}) aims to solve a well-posed problem. However, with $\Delta t < 0$,
together with {\em appropriate} data $U^0$ at time $T_{max}$, marching backward from $T_{max}$ in Eq.(\ref{eq:9.0006a}) attempts to solve an ill-posed problem.
Define the following norms for  three component vectors such as $W(.,t)$ and $U^n$,
\begin{equation}
\begin{array}{l}
\parallel W(\cdot,t) \parallel_2= \left\{\parallel u(\cdot,t) \parallel_2^2 + \parallel v(\cdot,t) \parallel_2^2 +\parallel w(\cdot,t)  \parallel_2^2\right\}^{1/2}, \\
\ \\
\parallel U^n \parallel_2= \left\{\parallel u^n \parallel_2^2 + \parallel v^n \parallel_2^2 + \parallel w^n \parallel_2^2 \right\}^{1/2}, \\
\ \\
||| W |||_{2,\infty}=\sup_{~0 \leq t \leq T_{max}} \left \{\parallel W(\cdot,t) \parallel_2 \right \}. \\
\end{array}
\label{eq:9.0006b}
\end{equation}
\newtheorem{lambdaQ}{Lemma}
\begin{lambdaQ}
With $p >1$, and $~\zeta_m, ~q_m,$ as in Eq.~(\ref{eq:9.0005a}), fix a positive integer $J,$
and choose $\omega \geq (\zeta_J)^{1-p}$. Then,
\begin{equation}
q_m \left(1+ |\Delta t| \zeta_m \right) \leq 1 + |\Delta t| \zeta_J, \qquad m\geq 1.
\label{eq:9.0005a1}
\end{equation}
\end{lambdaQ}
\ \\
{\em Proof~}: See  \cite[Lemma 6]{Thermo}.
\ \\
\newtheorem{newStable}[lambdaQ]{Lemma}
\begin{newStable}
With $\omega,~p,~\zeta_J,$ as in Lemma 1, and $R$ as in Eq.(\ref{eq:9.0006a}), we have $\parallel R \parallel_2 \leq 1 + |\Delta t| \zeta_J$. The
explicit scheme in Eq.(\ref{eq:9.0006a}) is unconditionally stable, and
\begin{equation}
\parallel U^{n} \parallel_2 =\parallel R^n U^0 \parallel_2 \leq \exp\{ n|\Delta t|\zeta_J \} \parallel U^0 \parallel_2, \qquad n=1, 2, \cdots, N.
\label{eq:2.003}
\end{equation}
\end{newStable}
\ \\
{\em Proof~}: See  \cite[Lemma 7]{Thermo}.
\ \\
\newtheorem{Prelim}[lambdaQ]{Lemma}
\begin{Prelim}
Let $W(t)$ be the unique solution of Eq.(\ref{eq:9.0004}).
Then, with $G,~S$ and $P$ as in Eq.(\ref{eq:9.0006}), the definitions of the norms in Eq.(\ref{eq:9.0006b}),
and $0 \leq n \leq N$,
\begin{eqnarray}
\parallel \tau^n \parallel_2 &\leq& 1/2(\Delta t)^2~ |||G^2W |||_{2,\infty}, \nonumber \\
\parallel W^n-SW^n \parallel_2 &\leq& \omega |\Delta t| ~||| PW |||_{2,\infty}, \nonumber \\
|\Delta t| \parallel GW^n-SGW^n \parallel_2 &\leq& \omega (\Delta t)^2~||| PGW |||_{2,\infty}.
\label{eq:2.006}
\end{eqnarray}
\end{Prelim}
\ \\
{\em Proof~}: See  \cite[Lemma 3]{Thermo}.
\ \\
\newtheorem{Error}{Theorem}
\begin{Error}
With $\Delta t > 0$, let $W^{n}$ be the unique solution of Eq.~(\ref{eq:9.0004}) at $t=n \Delta t$. Let $U^{n}$
be the corresponding solution of the forward explicit scheme in Eq.~(\ref{eq:9.0006a}), and let $p, ~\zeta_J, ~\omega,$ be as in Lemma 1.
If $ER(t)\equiv U^{n}-W^{n},$ denotes the error at $t= n \Delta t,~~n=0, 1, 2, \cdots, N,$ we have
\begin{eqnarray}
~~~~&~&~\parallel ER(t) \parallel_2 \leq e^{t \zeta_J} \parallel ER(0) \parallel_2
+\left \{\omega(e^{t \zeta_J}-1)/\zeta_J \right \} |||PW|||_{2,\infty} \nonumber \\
~~~~&+& \left \{(e^{t \zeta_J}-1)/\zeta_J \right \} \left \{\omega \Delta t~ |||PGW|||_{2,\infty} +(\Delta t/2)~ |||G^2W|||_{2,\infty} \right \}.
\label{eq:92.007}
\end{eqnarray}
\end{Error}
\ \\
{\em Proof~}: See  \cite[Theorem 1]{Thermo}.
\ \\
\ \\
In the above  forward problem, as $\Delta t \downarrow 0$, we are left with the error $e^{t \zeta_J} \parallel ER(0) \parallel_2$, originating in possibly erroneous initial data,
together with the {\em stabilization penalty}, represented by the second term in Eq, (\ref{eq:92.007}). These errors grow monotonically as $t$ increases. 
If $T_{max}$ is large, the accumulated distortion may become unacceptably large as $t\uparrow T_{max}$, and the stabilized explicit scheme is not useful in that case.

Marching backward from $t=T_{max}~$ in the backward problem,
solutions exist only for a restricted
class of data satisfying certain smoothness constraints. Such data are seldom known with sufficient accuracy. This is especially true of the hypothetical data $W^*(\cdot,T_{max})$ in the present data assimilation problem. It will be assumed that the given data
$U_b= W^*(\cdot,T_{max})$,  differ from the necessary exact data $W(\cdot,T_{max})$, by an {\em unknown} amount $\delta$ in the ${\cal{L}}^2(\Omega)$ norm.
\begin{equation}
\parallel U_b -W (\cdot,T_{max})\parallel_2 \leq \delta.
\label{eq:2.0101A}
\end{equation}
This leads to the following result, as proved in \cite[Theorem 2]{Thermo}.
\ \\
\newtheorem{ErrorB}[Error]{Theorem}
\begin{ErrorB}
With $\Delta t < 0$, let $W^{n}$ be the unique solution of the forward well-posed problem in Eq.~(\ref{eq:9.0004}) at $s=T_{max}-n|\Delta t|$. Let $U^{n}$
be the solution of the backward explicit scheme in Eq.~(\ref{eq:9.0006a}), with initial data $U(0)=U_b=[f_b,g_b,h_b]$ as in Eq.(\ref{eq:2.0101A}).
Let $p, ~\zeta_J, ~\omega,$ be as in Lemma 6.
If $ER(s)\equiv U^{n}-W^{n},$ denotes the error at $s=T_{max}- n |\Delta t|,~~n=0, 1, 2, \cdots, N,$ then
\begin{eqnarray}
~~~~&~&~\parallel ER(s) \parallel_2 \leq  \delta e^{n|\Delta t| \zeta_J} + \left \{\omega (e^{n|\Delta t| \zeta_J}-1)/\zeta_J \right \} |||PW|||_{2,\infty} \nonumber \\
~~~~&+& \left \{(e^{n|\Delta t| \zeta_J}-1)/\zeta_J \right \}\left\{\omega |\Delta t|~ |||PGW|||_{2,\infty} +(|\Delta t|/2)~ |||G^2W|||_{2,\infty} \right \}.
\label{eq:92.0101} 
\end{eqnarray}
\end{ErrorB}
\ \\
\subsection{Application to data assimilation}
In Theorems 1 and 2 above, define the constants $K_1$ through $K_5$ as follows, and consider the values shown in Table 1 below.
\begin{eqnarray}
        ~~~&~&~K_1=e^{\zeta_J T_{max}},~~K_2=\omega (e^{\zeta_J T_{max}}-1)/\zeta_J,~~K3=|\Delta t|K_2,~~K_4=K_3/(2 \omega), \nonumber \\
        ~~~&~&~K_5=K_2 |||PW|||_{2,\infty} + K_3 |||PGW|||_{2,\infty} +K_4 |||G^2W|||_{2,\infty}.
\label{eq:3.055}
\end{eqnarray}
\begin {center}
\ \\
{\em TABLE 1} \\
{\em Values of $K_1$ through $K_4$ in  Eq.~(\ref{eq:3.055}), with  following  parameter values:}\\
{$T_{max}=1.6 \times10^{-4},~|\Delta t|=(4/3)\times10^{-7},~p=3.35,~\zeta_J=19800,~\omega=\zeta_J^{(1-p)}= 8\times10^{-11}$.} \\
\ \\
\begin{tabular}{|c|c|c|c|} \hline
        {$K_1= e^{\zeta_J T_{max}}$} & {$K_2=(\zeta_J)^{-p} (K_1-1)$} & {$K_3=|\Delta t|K_2$}& {$K_4=K_3/(2 \omega)$} \\ \hline
	{$K_1 < 23.8$} & {$K_2 < 9.2  \times 10^{-14}$} & {$K_3 < 1.3 \times 10^{-20}$}&{$K_4 < 7.7 \times 10^{-11}$} \\ \hline
\end{tabular}
\end{center}
\ \\
\ \\
As outlined in the Introduction, data assimilation applied to the system in Eq.~(\ref{eq:9.0004}), is the problem of finding initial values $[u(.,0), v(.,0), w(.,0)]$, at $t=0$,
that can evolve into useful approximations to $W^*(\cdot,T_{max})$, the given hypothetical
data at an appropriate time $T_{max} > 0$. If the true solution in  Eq.~(\ref{eq:9.0004}) does not have exceedingly large values for $|||PW|||_{2,\infty},~ |||PGW|||_{2,\infty}$,~ or $|||G^2W|||_{2,\infty}$,
the parameter values chosen in Table 1, together with Theorem 2, indicate that marching backward to time $t=0$ from the hypothetical data $W^*$ at $T_{max}$, leads to an error $ER(0)$, satisfying
\begin{equation}
        \parallel ER(0) \parallel_2 \leq \delta K_1 + K_5,
        \label{eq:3.056}
\end{equation}
where the constant $K_5$ may be negligibe compared to $\delta K_1$. Next, from Theorem 1, marching forward to time $T_{max}$ using the inexact computed initial values $U(\cdot,0)$, leads to an error $ER(T_{max})$, satisfying
\begin{equation}
        \parallel ER(T_{max}) \parallel_2 \leq K_1(\delta K_1 + K_5) + K_5.
        \label{eq:3.057}
\end{equation}
The error $ ER(T_{max})$ in Theorem 1 is the difference at time $T_{max}$, between the unknown unique solution $W(\cdot, t)$ in Eq.(\ref{eq:9.0004}), and the computed numerical approximation to it, $U(\cdot,t)$,  provided by the stabilized forward explicit scheme.
However, $\parallel W^*(\cdot,T_{max}) - W(\cdot,(T_{max}) \parallel_2 \leq \delta$, if $W^*(\cdot,T_{max})$ is the given hypothetical data. Hence, using
the triangle inequality, one finds
\begin{equation}
        \parallel W^*(\cdot,T_{max}) - U(\cdot,T_{max}) \parallel_2 \leq \delta (1+K_1^2) +K_5 (1+K_1).  
        \label{eq:3.058}
\end{equation}
Therefore, data assimilation is successful only if the inexact computed initial values $U(\cdot,0)$ at $t=0$, lead to a sufficiently small right hand side
in  Eq.(\ref{eq:3.058}). Clearly, the value of $\zeta_J T_{max}$, together with the unknown value of $\delta$, will  play a vital role.
From Table 1, we find $\delta (1+K_1^2) <  568 ~\delta$. However, with $T_{max}$ chosen five times larger, one would find 
$\delta (1+K_1^2) > (5.8 \times 10^{13}) ~\delta$.
\ \\

\subsection{Using the Laplacian for smoothing when $L$ has variable coefficients} All of the discussion and results in \cite[Section 6]{Thermo} on using the Laplacian for smoothing, can be applied to the linear system in Eq.~(\ref{eq:9.0004}). With $\rho,~\Lambda,~Q$ as in Eq.~(\ref{eq:9.0005a}),
let $\Gamma=\rho(I-\Delta)$. For real $q > 1$, and $\epsilon > 0,$ define $~Q_{\Delta}=\exp\{-\epsilon |\Delta t| \Gamma^q\}$. In domains where closed form expressions for the eigenfunctions of the Laplacian are known, it may be 
advantageous to use the smoothing operator $Q_{\Delta}$ in lieu of $Q$ in the
stabilized explicit scheme in Eq.~(\ref{eq:9.0006a}). This is feasible for
differential operators $L$ where the hypothesis in \cite[Eq. (6.2)]{Thermo} 
is valid, so that, with appropriately chosen $(\epsilon, q)$, one has $\parallel Q_{\Delta}g \parallel_2 \leq \parallel Q g \parallel_2$, for all $g \in {\cal{L}}^2(\Omega)$ and sufficiently small $|\Delta t|.$ When this is the case, 
Theorems 1 and 2 remain valid, with $S_{\Delta}$ and $P_{\Delta}$ replacing $S$ and $P$.
Moreover, as described in \cite[Section 6.1]{Thermo}, and as will be demonstrated in the 
computational experiments discussed below, 
it may be possible to use efficient FFT algorithms to synthesize $Q_{\Delta}$, even in problems defined on non-rectangular domains $\Omega$.

\section{Data assimilation in nonlinear coupled sound and heat flow in non-rectangular region, using FFT Laplacian smoothing}
%\begin{figure}
%	\centerline{\includegraphics[width=5.25in]{LincolnIntensitySurface.eps}}
%	\caption{Nonlinear coupled sound and heat flow experiment in Section 10 lies outside scope of linear theory }
%\end{figure}
\begin{figure}
\centerline{\includegraphics[width=5.25in]{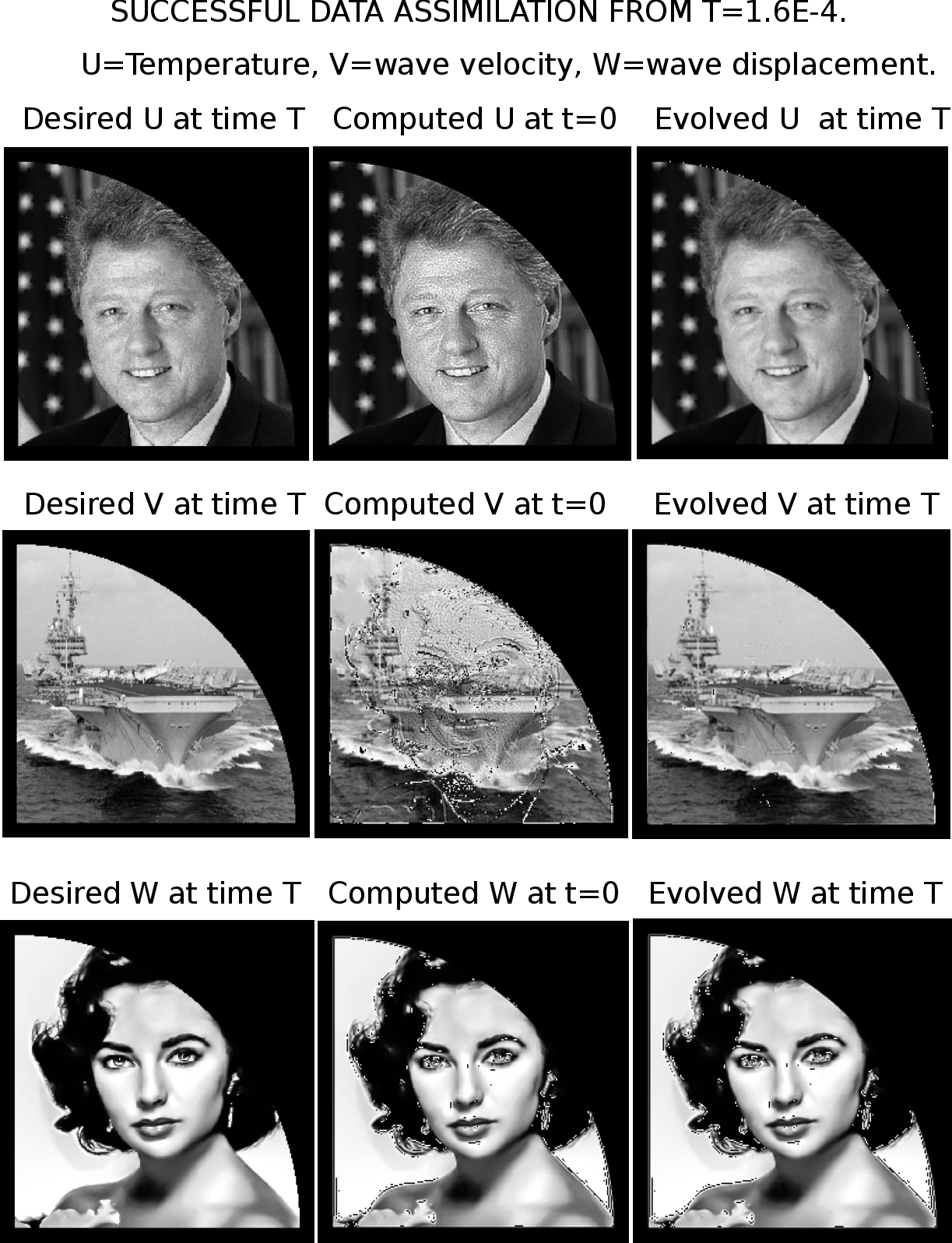}}
\caption{FIGURE 2. Successful data assimilation experiment. See summary in Table 2. Above nonlinear coupled sound and heat flow experiment lies outside scope of linear theory developed in Section 2. As explained in in the discussion following Eq.~(\ref{eq:309.0001}), enclosing quarter circle region $\Omega$ in unit square
	$\Psi$, allows use of FFT Laplacian smoothing operator $Q_{\Delta}$, in backward reconstruction with scheme in Eq.~(\ref{eq:9.0006a}).}
\end{figure}
\begin{figure}
\centerline{\includegraphics[width=5.25in]{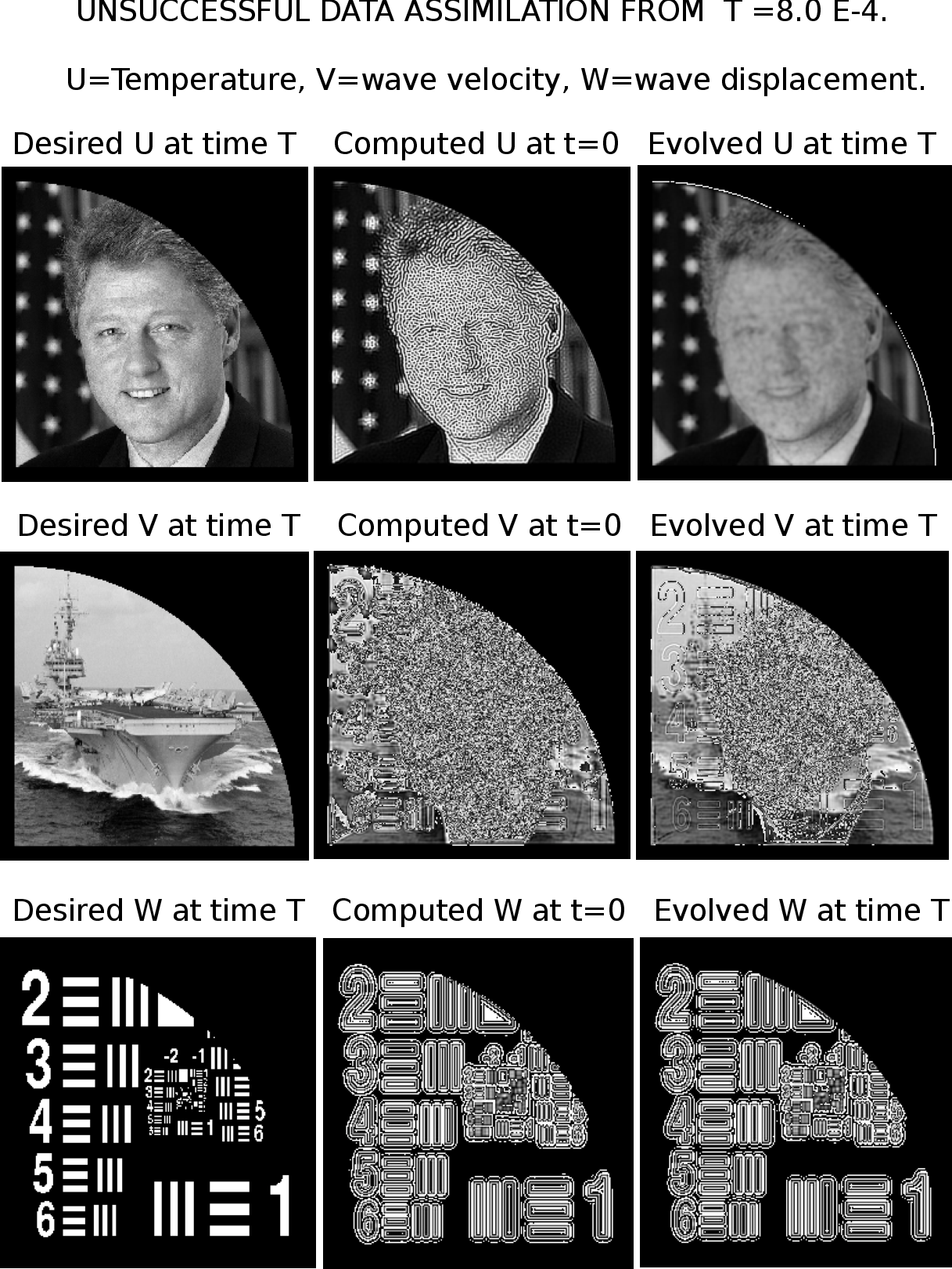}}
	\caption{FIGURE 3. Failure of data assimilation with significantly larger value for $T_{max}$. See summary in Table 3. As previously explained, enclosing quarter circle region $\Omega$ in unit square
	$\Psi$, allows use of FFT Laplacian smoothing operator $Q_{\Delta}$, in backward reconstruction with scheme in Eq.~(\ref{eq:9.0006a}).}
\end{figure}

We now highlight the versatility of the stabilized  scheme in Eq.~(\ref{eq:9.0006a}), by considering a nonlinear example in a non-rectangular region, in which the explicit computation proceeds by lagging the nonlinearity at the previous
time step.
Let $\Omega$ be the open quarter circle region in the $(x,y)$ plane,
\begin{equation}
0.05 < x, y < 0.95,~~(x-0.05)^2 + (y-0.05)^2 < (0.9)^2,
\label{eq:98.0001}
\end{equation}
let $T_{max}=1.6 \times 10^{-4}$, and let  $L$ be the nonlinear differential operator defined as follows on functions $z(x,y,t)$ on $\Omega \times (0,T_{max})$:
\begin{equation}
Lz=-0.00085~s(z)~\nabla. \{q(x,y) \nabla z\} -  2.75 (z_x + z_y), \\
\label{eq:99.00001X}
\end{equation}
where
\begin{equation}
s(z)=\exp\{0.005 z \}, ~~~~~~1 < q(x,y)=\left\{1 + 2 \sin \pi x \sin \pi y \right\} \leq 3,
\label{eq:99.000078X}
\end{equation}
With $a = 6,~b=5,~d=0.95,$ and $(x,y,t) \in \Omega \times (0, T_{max})$, consider the system
\begin{equation}
\begin{array}{l}
u_t=-b L u - dv, \\
\ \\
v_t = aL u -a L w, \\
\ \\
w_t= v, \\
\ \\
	u(x,y,T_{max})=f(x,y),~~v(x,y,T_{max})=g(x,y),~~w(x,y,T_{max})=h(x,y), \\
\ \\
u(x,y,t)=v(x,y,t)=w(x,y,t)=0, ~~(x,y,t) \in \partial \Omega \times [0, T_{max}].
\end{array}
\label{eq:309.0001}
\end{equation}

In Eq.~(\ref{eq:309.0001}), $u(x,y,t)$ denotes the temperature, $w(x,y,t)$ is the wave displacement, and $v=w_t(x,y,t)$ is the wave velocity. The hypothetical data at time $T_{max}$, namely, $f(x,y),~ h(x,y),~ g(x,y),$ are the three
images shown in the leftmost column in Figure 2.  Here, the quarter circle region $\Omega$ was enclosed in the unit square $\Psi=\{0 < x,y < 1\}$. A $512 \times 512$ uniform grid was
imposed on $\Psi$, leading to a discrete boundary $\partial\Omega_d$ consisting of the grid points closest to $\partial\Omega$. This was assumed to sufficiently well-approximate  $\partial \Omega$. With $\Delta x= \Delta y= 1/512$, and
$\Delta t = -(4/3) \times 10^{-7}$, homogeneous boundary conditions were applied on $\partial\Omega_d$. Lagging the nonlinearity at the previous time step, explicit time differencing together with centered finite differencing in the space variables,
were used in the {\em stabilized backward computation}, as described below, for $1200~\Delta t$. This produced the images at $t=0$, shown in the middle column in Figure 2. 
The actual computed data at $t=0$ involve negative values. These values are not used in forming and displaying the middle column images, but are nevertheless retained as necessary to enable computation of
the images in the rightmost column in Figure 2.

With $\rho$ as in Eq.~(\ref{eq:9.0005a}),
$\Gamma=\rho(I-\Delta)$, real $q > 1$, and $\epsilon > 0,$ let $~Q_{\Delta}=\exp\{-\epsilon |\Delta t| \Gamma^q\}$.
In the above stabilized backward computation 
at each time step $m$ in Eq.~(\ref{eq:9.0006a}), after applying the operator $(I + \Delta t G)$ to $U^m$ on $\Omega \subset \Psi$, the solution is extended to all of $\Psi$ by defining it to be zero on $\Psi-\Omega$.
FFT algorithms are then applied on $\Psi$ to synthesize $Q_{\Delta}$, and produce $U^{m+1}=S_{\Delta}(I+\Delta t G)U^m$, while retaining only the values of $U^{m+1}$ on $\Omega$.  This process is then repeated at the next time step. Interactive trials are needed to locate appropriate values for $(\epsilon, q)$. Here, a parameter pair $\epsilon= 8.0 \times 10^{-11},~q=3.35$, was arrived at after very few trials.

As is evident from Figure 2, together with the accompanying data in Table 2, data assimilation was successful in that computational experiment. At time $T_{max}$, the evolved $L^1$ norms closely match the 
desired $L^1$ norms, and the resulting $L^1$ relative errors are small. In the linear selfadjoint problem contemplated in Table 1 and  Eq.(\ref{eq:3.058}), we found  $\delta (1+K_1^2) < 568~ \delta$, when using
the same parameter values that were used in the nonlinear problem in Figure 2. However, with $T_{max}$ chosen five times larger in the linear problem, we found $\delta (1+K_1^2) >(5.8 \times 10^{13})~ \delta$.

We now consider data assimilation for the same nonlinear problem described in Eq.~(\ref{eq:309.0001}), with the same parameters used in Figure 2 and Table 2, except for $T_{max}$, which is now chosen 
{\em five time larger}, with the Elizabeth Taylor image replaced by the USAF 1951 resolution chart. As can be seen from Figure 3 and Table 3, data assimilation is now unsuccessful. 
The images in the rightmost column in Figure 3, do not match the desired images in the leftmost column.
The evolved $L^1$ norms are not
good approximations to the desired $L^1$ norms at time $T_{max}$, and the accompanying $L^1$ relative errors are quite large. While the failure in the above nonlinear
experiment with $T_{max}=8.0 \times 10^{-4}$ is less catastrophic than was predicted in the linear selfadjoint case with the same $T_{max}$ value,
the linear analysis in Table 1 and  Eq.(\ref{eq:3.058}), was a useful guide. As previously noted, 
successful reconstruction in backward dissipative evolution equations is necessarily limited, given the associated uncertainty estimates, in \cite{karen, payne3, knopslog, knops,payne,car00,car77,hao}.

\begin {center}
\ \\
{\em TABLE 2} \\
{\em Behavior of $L^1$ norm in successful data assimilation example shown in Figure 2.}\\
{$T_{max}=1.6 \times10^{-4},~|\Delta t|=(4/3)\times10^{-7},~p=3.35,~\zeta_J=19800,~\omega=\zeta_J^{(1-p)}= 8\times10^{-11}$.} \\
\ \\
\begin{tabular}{|c|c|c|c|} \hline
	{Image~~~~} & {Desired $L^1$ norm }& {Evolved $L^1$ norm} & {$L^1$ relative error} \\ \hline
        {Bill Clinton} & {61.06~($T_{max}$)} & {60.80~($T_{max}$)}&{7.46~\%~($T_{max}$)} \\ \hline
        {USS Kitty Hawk} & {93.60~($T_{max}$)} & {93.07~($T_{max}$)}&{4.41 \%~($T_{max}$)} \\ \hline
        {Elizabeth Taylor} & {99.95~($T_{max}$)} & {99.27~($T_{max}$)}&{2.69 \%~($T_{max}$)} \\ \hline
\end{tabular}
\end{center}
\ \\
\ \\
\begin {center}
\ \\
{\em TABLE 3} \\
{\em Behavior of $L^1$ norm in unsuccessful data assimilation example shown in Figure 3.}\\
{$T_{max}=8.0 \times10^{-4},~|\Delta t|=(4/3)\times10^{-7},~p=3.35,~\zeta_J=19800,~\omega=\zeta_J^{(1-p)}= 8\times10^{-11}$.} \\
\ \\
\begin{tabular}{|c|c|c|c|} \hline
	{Image~~~~} & {Desired $L^1$ norm }& {Evolved $L^1$ norm} & {$L^1$ relative error} \\ \hline
	{Bill Clinton} & {61.06 ~($T_{max}$)} & {56.57~($T_{max}$)}&{18.11 \%~($T_{max}$)} \\ \hline
	{USS Kitty Hawk} & {93.60~($T_{max}$)} & {85.81~($T_{max}$)}&{48.94 \%~($T_{max}$)} \\ \hline
	{USAF chart} & {39.27~($T_{max}$)} & {41.59~($T_{max}$)}&{26.25 \%~($T_{max}$)} \\ \hline

\end{tabular}
\end{center}
\ \\
\ \\

\section{Concluding Remarks}
With proper parameter choices, stabilized explicit schemes appear to be helpful in difficult data assimilation problems, involving
non differentiable data and nonlinear dissipative systems. 
Examples of failure in such problems are also instructive and valuable.

Along with \cite{Burg,Advec,Thermoassim, LeapNavAssim}, the results in the present paper invite  useful scientific debate and comparisons,
as to whether equally good or better
results might be achieved, using the computational methods described in \cite{cintra, howard, blum, qizhi, arcucci, Antil, Chong, Lund,auroux1,ou,auroux2,auroux3,pozo,gosse,gomez, xu, tomislava,camposvelho}. As an alternative computational approach, backward marching stabilized explicit schemes may also be helpful, if needed,
in validating computations involving  machine learning.

\ \\

\end{document}